\documentclass[twoside,11pt]{amsart}
\usepackage{amsmath,latexsym,amssymb,times,enumerate}

\usepackage{verbatim}

\def\demo{\noindent{\bf Proof. }}
\def\sqr#1#2{{\vcenter{\hrule height.#2pt
        \hbox{\vrule width.#2pt height#1pt \kern#1pt
                \vrule width.#2pt}
        \hrule height.#2pt}}}
\def\square{\mathchoice\sqr64\sqr64\sqr{4}3\sqr{3}3}
\def\QED{\hfill$\square$}

\def\lto{\longrightarrow}

\def\NN{\mathbb N}

\def\a{{\mathfrak a}}

\def\m{{\mathfrak m}}

\def\p{{\mathfrak p}}

\def\B{\mathcal B}
\def\A{\mathcal A}
\def\H{\mathcal H}

\def\C{\mathcal C}
\def\J{\mathcal J}

\def\R{\mathcal R}
\def\F{\mathcal F}
\def\D{\mathcal D}
\def\E{\mathcal E}
\def\ux{\underline{x}}
\def\uj{\underline{j}}
\def\ui{\underline{i}}

\def\o1{{\overline 1}}

\def\ux{{\underline{x}}}
\def\uk{{\underline{k}}}
\def\ut{{\underline{t}}}

\def\ms{\medskip}
\def\bs{\bigskip}
\def\s{\smallskip}

\def\gb{Gr\"{o}bner \,}

\newtheorem{Theorem}{Theorem}[section]
\newtheorem{Lemma}[Theorem]{Lemma}
\newtheorem{Corollary}[Theorem]{Corollary}
\newtheorem{Proposition}[Theorem]{Proposition}

\newtheorem{Remark}[Theorem]{Remark}
\newtheorem{Definition}[Theorem]{Definition}

\newtheorem{Notation and Discussion}[Theorem]{Notation and Discussion}
\newtheorem{Assumptions and Discussion}[Theorem]{Assumptions and Discussion}
\newtheorem{Assumptions}[Theorem]{Assumptions}
\newtheorem{Example}[Theorem]{Example}

\newtheorem{Strategy}[Theorem]{Strategy}

\begin{document}

\baselineskip=16pt


\title[Rees Algebras of Truncations of Complete Intersections]
{\Large\bf Rees Algebras of Truncations of Complete Intersections}
\author[K. N. Lin \and C. Polini]
{Kuei-Nuan Lin  \and Claudia Polini}

\thanks{AMS 2010 {\em Mathematics Subject Classification}.
Primary 13A30; Secondary 13H15, 13B22, 13C14, 13C15, 13C40.}
\thanks{The first author was partially supported by an AWM-NSF mentor travel grant.}

\thanks{ The second author was partially supported by NSF grant DMS-1202685 and NSA grant H98230-12-1-0242.
}


\address{
Department of Mathematics and Statistics, 
Smith College, Northampton, Massachusetts 01063
}\email{klin@smith.edu}
\address{Department of Mathematics, University of Notre Dame, Notre Dame, Indiana 46556} \email{cpolini@nd.edu}

\vspace{-0.1in}

\begin{abstract}  In this paper we describe the defining equations of the Rees algebra and the special fiber ring of a truncation $I$ of a complete intersection ideal in a polynomial ring over a field with homogeneous maximal ideal ${\mathfrak m}$.  To describe explicitly the Rees algebra $\mathcal R(I)$ in terms of generators and relations we  map another Rees ring $\mathcal R(M)$ onto it, where $M$ is  the direct sum of powers of $\m$. We compute a  Gr\"{o}bner basis of the ideal defining $\mathcal R(M)$. It turns out that the normal domain  $\mathcal R(M)$ is a Koszul algebra and from this we deduce that  in many instances $\mathcal R(I)$ is a  Koszul algebra as well. \end{abstract}

\maketitle

\vspace{-0.2in}

\section{Introduction}

\s
In this paper we investigate the Rees algebra $\mathcal R(I)= R[It]$ as well as the special fiber ring $\mathcal F(I)=\mathcal  R(I) \otimes k$ of an ideal $I$ in a standard graded algebra $R$ over a field $k$. These objects are important to commutative algebraists  because they encode the asymptotic behavior of the ideal $I$ and to algebraic  geometers because  their projective schemes define the blowup and the special fiber of the blowup of the scheme ${\rm
Spec}(R)$ along  $V(I)$. 
One of the central problems in the theory of Rees rings is to describe $\mathcal R(I)$ and $\mathcal F(I)$ in terms of generators and relations (see for instance \cite{RV, V, GG, UV, SUV0, M, MU, J, H, CHW, HSV, HSV2}). This is a challenging quest which is open for most classes of ideals, even three generated ideals in a polynomial ring in two variables (see for instance \cite{B, BD, KPU3}). The goal is  to find an ideal $\A$ in a polynomial ring $S=R[T_1, \ldots,T_s]$ so
that $\mathcal R(I) =S/\A$.    

If the ideal $I$ is generated by forms of the same degree, then these forms define rational maps between projective spaces and the special fiber ring and the Rees ring describe the image and the graph of such rational maps, respectively.  By computing the defining equations of these algebras, one is able to exhibit the implicit equations of
the graph and  of the variety parametrized by the map. This 
classical and difficult problem in elimination theory has also been studied
in applied mathematics, most notably in modeling theory, where it is known as the implicitization problem (see for instance \cite{BuC, BuJ, CS,EM}).

If the ideal $I$ is not generated by forms of the same degree,  one can consider  the truncation of $I$ past its generator degree. In this paper we treat truncations of complete intersection ideals in a polynomial ring. More precisely, let $R=k[x_1,\ldots, x_n]$ be a polynomial ring over a field $k$ with homogeneous maximal ideal ${\mathfrak m}$, 
let $f_{1}$, \ldots, $f_{r}$ be a homogeneous regular sequence in $R$ of degrees $d_{1}\ge \ldots \ge d_{r}$, let $d\ge d_1$ be an integer, and write  $a_i=d-d_i$.  The truncation $I=(f_1, \ldots, f_r)_{\ge d}$ of the complete intersection $(f_1,\ldots, f_r)$ in degree $d$ is  the $R$-ideal generated by the forms in $(f_1,\ldots, f_r)$ of degree at least $d$. In other words, $I=(f_1, \ldots, f_r) \cap \m^d= \sum_{i=1}^r \mathfrak{m}^{a_{i}}f_{i} $. The Cohen-Macaulayness of the Rees algebra of such ideals was previously studied in \cite{HT}, where the authors show that $\mathcal R(I)$ is Cohen-Macaulay for all $d>D$ and they give a sharp estimate for $D$. However, the defining equations of $\R(I)$ were unknown. In this paper we describe them explicitly for all $d$ when $r=2$ and for $d\ge d_1+d_2$ when $r\ge3$. 
Furthermore we prove that the Rees ring and the special fiber ring are Koszul algebras  for $d\ge d_1+d_2$ and for  $d\ge d_1+d_2-1$ if $r=2$.

To determine the defining equations of $\R(I)$ we  map another Rees ring $\mathcal R(M)$ onto $\mathcal R (I)$, 
where $M$ is the module $\m^{a_1} \oplus \ldots \oplus\m^{a_r}$. Our aim then becomes to find the defining ideal of $\mathcal R(M)$ and the kernel $Q$, 
\[ 0 \rightarrow Q \longrightarrow \mathcal R(M) \longrightarrow \mathcal R(I) \rightarrow 0\, .
\]

The problem of computing the implicit equation of  $\mathcal R (M)$ is interesting in its own right and it was previously addressed in \cite{R}, where the  relation type of $\mathcal R(M)$ was computed. We solve it in Section~\ref{section1}. It  turns
out that $\mathcal R(M)$ and $\mathcal F(M)$ are normal domains whose defining ideals have a Gr\"{o}bner basis of quadrics; hence, they are Koszul algebras. For $r=2$, the kernel $Q$ is a height one prime ideal of the normal domain $\mathcal R(M)$; therefore it is  a divisorial ideal of $\mathcal R(M)$. Our goal is then reduced to  explicitly describing ideals that represent the
elements in the divisor class group of $\mathcal R(M)$. The approach is very much inspired by \cite{KPU2} and \cite{KPU1}. For $r=2$ and $d \ge d_1+d_2-1$ or $r\ge 3$ and $d \ge d_1 +d_2$, the ideal $I$ has a linear presentation and $Q$ turns out to be a linear ideal in the $T's$. Using this we prove  that the defining ideal of $\mathcal R(I)$ has a quadratic Gr\"{o}bner basis and hence  $\mathcal R(I)$ is a Koszul algebra as well.

\bs

\section{The Blowup algebras of direct sums of powers of the maximal ideal}\label{section1}

Let $R=k[x_1,\ldots, x_n]$ be a polynomial ring over a field $k$ with homogeneous maximal ideal ${\mathfrak m}$.  
Let $0 \le a_{1} \le \ldots \le a_{r}$ be integers. Write $\a=a_1, \ldots, a_r$.   In this section, we will describe explicitly  the 
Rees algebra and the special fiber ring of the module $M=M_{\a}=\mathfrak{m}^{a_{1}}\oplus\cdots\oplus\mathfrak{m}^{a_{r}}$ 
in terms of generators and relations and we will prove that they are Koszul normal domains. We will end the section with a study of the divisor class  group of the blowup algebras of $M$.

\begin{Definition}\label{N0}
{\rm Let $R=k[x_1,\ldots, x_n]$ be a polynomial ring over a field $k$.  
Let $a$ be a positive integer, write $\J_a$  and $\J'_a$ for the two sets of multi-indices in  $(\NN \cup\{0\})^{n-1}$ and $\NN^{n-1}$, respectively, that are defined as follows
\[\J_a=\{ \uj=(j_{n-1}, \ldots, j_1) \ | \  0 \le j_1\le \ldots \le j_{n-1}\le a\}\, , \]
\[\J_a'=\{ \uj=(j_{n-1}, \ldots, j_1) \ | \  1 \le j_1\le \ldots \le j_{n-1}\le a\}\, . \]
Write $\ux^{a,\uj} $  and $\ux^{a,\uj,s} $ for the monomials \[\hspace{.5cm}\ux^{a,\uj} =\prod_{i=1}^{n} x_i^{j_i-j_{i-1}}\] with $ \uj\in \J_a, j_0=0,  j_n=a$  and  \[\ux^{a,\uj,s} =\frac{x_s}{x_1}  \ux^{a,\uj} \]  with $ \uj\in \J_a',  1\le s \le n \, .$
}
\end{Definition}
The Rees algebra $ \mathcal{R}(M)$ of $M$  is the subalgebra
\[
\hspace{1cm} \mathcal{R}(M)=R [\{\ux^{a_l, \uj}t_l \mid 1 \le l \le r, \uj\in \J_{{a_l}}\} ] \subset R[t_1, \ldots, t_r]
\]
of the polynomial ring $R[t_1, \ldots, t_r]$, while the special fiber ring $\F(M)$ is the subalgebra of $\R(M)$ 
\[
\F(M)= k[ \{ \ux^{a_l, \uj}t_l \mid 1 \le l \le r, \uj\in \J_{{a_l}} ] \subset \R(M)\, .
\] 
To find a presentation of these algebras we consider the polynomial rings 
\begin{eqnarray*}
T=T_{\a}&=&k[\{T_{{l},\underline{j}} \mid 1 \le l \le r, \  \underline{j} \in \J_{{a_l}}\}] \\ 
S=S_{\a}&=&R \otimes_k T_{\a} =T_{\a}[x_1, \ldots, x_n]
\end{eqnarray*}
in the new variables $T_{{l}, \underline{j}}$ and 
 the epimorphisms of algebras
\begin{equation}\label{map}
\phi:S\rightarrow\mathcal{R}(M) \qquad \psi: T\rightarrow \F(M)
\end{equation}
defined by \[\phi(T_{{l}, \underline{j}})=\psi(T_{{l}, \underline{j}})=\ux^{a_l, \uj}t_l\, .\] Notice that $\psi$ is the restriction of $\phi$ to $T$. 

We can assume that the $a_i$'s are all positive because if $a_{i}=0$ for $1\le i \le s$ with $s \le r$ then  the Rees algebra $\R(M)$ is isomorphic to a polynomial ring over the Rees algebra of $\mathfrak{m}^{a_{s+1}}\oplus\cdots\oplus\mathfrak{m}^{a_{r}}$
\[\mathcal{R}(\oplus_{l=1}^{r}\mathfrak{m}^{a_{l}})\cong\mathcal{R}(\oplus_{l=s+1}^{r}\mathfrak{m}^{a_{l}})[t_{1},\cdots,t_{s}] \, ,\]
and likewise for the special fiber ring. 
Furthermore, we can treat simultaneously the special fiber ring and the Rees algebra since \[\mathcal{R}(\oplus_{l=1}^{r}\mathfrak{m}^{a_{l}})\cong\mathcal{F}(\m \oplus (\oplus_{l=1}^{r}\mathfrak{m}^{a_{l}}))  \quad \text{and} \quad  \mathcal{F}(M)=\mathcal {R}(M)/\m\mathcal{R}(M)\, .
\]

\begin{Definition}\label{Def1}{\rm  Let $\tau$ be  the  lexicographic order on a set of  multi-indices in $(\NN \cup\{0\})^{n}$, i.e. 
$\underline{p}>\underline{q}$ if  the first nonzero entry of  $\underline{p}- \underline{q}$
is positive. Let $1 \le a_{1} \le \ldots \le a_{r}$ be integers.  Write $\a=a_1, \ldots, a_r$.
Order the set of multi-indices
$\{ (l, \uj)  \ | \ 1 \le l \le r , \uj \in \J'_{{a_l}}
\}$ by $\tau$. Write $T_{{l}, \uj, s}$ for the variable $T_{{l},j_{n-1},\ldots, j_s, j_{s-1}-1, \ldots, j_{1}-1}$. 



\begin{enumerate}
\item Let  $B_{\a}$ be the $n \times (\sum_{l=1}^{r}\binom{a_{l}+n-2}{n-1})$ matrix whose entry in the $s$-row and the $(l, \uj)$-column is the variable $T_{{l}, \uj, s}$ with $1\le s \le n$, $1 \le l \le r$, and $\uj \in \J'_{{a_l}}$.


\item  Let  $C_{\a}$ be the $n \times (1+\sum_{l=1}^{r}\binom{a_{l}+n-2}{n-1})$ matrix
\[
C_{\a}=\left[\begin{array}{lcl}
x_1 & \vline & \\
\vdots &\vline & B_{\a} \\
x_n & \vline &
\end{array}
\right] \, .
\]
\end{enumerate}
}
\end{Definition}

\begin{Example} {\rm For instance if $n=3$, $r=2$ and $\a=1,2$ then $B_{\a}$ is the $3 \times 4$ matrix 
\[B_{\a}=\left[\begin{array}{cccc}
T_{1,1,1} & T_{2,1,1} & T_{2,2,1} & T_{2,2,2} \\
T_{1,1,0} & T_{2,1,0} & T_{2,2,0} & T_{2,2,1} \\
T_{1,0,0} & T_{2,0,0} & T_{2,1,0} & T_{2,1,1} 
\end{array}
\right] \, .
\]
}
\end{Example}

\bigskip

Write  $\mathcal{L}$ for the kernel of the epimorphism $\phi$ defined in (\ref{map}). 
Notice that $\phi(T_{{l}, \uj, s})=\ux^{a,\uj,s} t_{l}=\frac{x_s}{x_1}  \ux^{a,\uj}t_{l}$
hence one can easily deduce the inclusion  
\[I_2(C_{\a}) \subset \mathcal{L} \, .\] 
Our goal is to show that the above inclusion is  an equality. In order to establish this claim  it suffices to show that  $I_2(C_{\a})$ is a prime ideal of dimension  at most $n+r ={\rm dim} \, \mathcal{R}(M)$ (see  for instance \cite[2.2]{SUV}). 

\begin{Theorem}\label{TH1} Let $R=k[x_1,\ldots, x_n]$ be a polynomial ring over a field $k$ with homogeneous  maximal ideal ${\mathfrak m}$. Let $1 \le a_{1} \le \ldots \le a_{r}$ be integers and let $M=\mathfrak{m}^{a_{1}}\oplus\cdots\oplus\mathfrak{m}^{a_{r}}$. The Rees algebra and the special fiber ring of $M$ are Koszul normal domains.  Furthermore, \[  \R(M)=S/I_2(C_{\a})\qquad \F(M)=T/I_2(B_{\a})\,  , \] where $\a$, $C_{\a}$, and $B_{\a}$ are as in Definition~\ref{Def1}.
\end{Theorem}

\ms

An important step in the proof of Theorem~\ref{TH1} is to show that the set of $2$ by $2$ minors of $C_{\a}$ forms a \gb basis for $I_2(C_{\a})$. In the next lemma we will show much more. Indeed, the set of $2\times2$ minors of any  submatrix $D_{\a}\, $ of $C_{\a}$ forms
a \gb basis for $I_{2}(D_{\mathfrak{a}})$. 
Notice that $\tau$ induces an ordering of the variables of $T$.   With respect to this ordering we consider the reverse lexicographic order  on the monomials in the ring $T$, which we  also call $\tau$. 
\ms

\begin{Remark}\label{C=B} {\rm Let $\a'=1, a_1, \ldots, a_r$ and denote with $T_{\a'}$ the polynomial ring associated with the sequence $\a'$. Notice that 
$T_{\a'}/I_2(B_{\a'})\cong S_{\a}/I_2(C_{\a})$ and the matrix 
 $C_{\a}$ is equal  (after changing the name of the variables)  to the matrix $B_{\a'}$.}
\end{Remark}
\ms

\begin{Lemma}\label{GbB} Adopt assumption ~\ref{Def1} and let $D_{\a}$ be any submatrix of $C_{\a}$ with the same number of rows.    The set of $2\times2$ minors of $D_{\a}\, $ forms
a \gb basis for $I_{2}(D_{\mathfrak{a}})$ with respect to $\tau$. 
 \end{Lemma} 

\ms

We use a general strategy to compute the \gb bases that we outline below.

\ms
\begin{Strategy}\label{strategy} {\rm Let $R$ be a polynomial ring  over a field with a fixed monomial order  $">"$. Let $ \mathcal{D}$ be a collection of $2\times2$ minors of a matrix with entries in $R$.
The Buchberger criterion  (see for instance \cite[15.8]{E}) states that a generating set $ \D$ of an ideal is a \gb basis for the ideal if the $S$-polynomials (or $S$-pairs) of any two elements of $ \D$ reduces to zero module $ \D$. The following strategy will be used to show that the remainders
of the $S$-pairs of the elements of $ \mathcal{D}$ reduces to zero modulo $ \D$ (see for instance  \cite[Section 2.9, Definition 1 and Theorem 3]{CLO}). Let
$ad-bc$ and $hg-ef$ be two elements of $ \mathcal{D}$. To compute the $S$-polynomial $S(ad-bc,hg-ef)$ we can assume, for instance, that $h=a$ and  $ad>bc$ and $ag>ef$ because otherwise the leading
terms of the two polynomials are relatively prime and therefore $S(ad-bc,ag-ef)$
reduces to zero modulo $ \mathcal{D}$ (see for instance \cite[Section 2.9, Proposition 4]{CLO}). 
We consider the matrix
\[
M=\left[\begin{array}{ccc}
d & b & 0\\
c & a & f\\
G & e & g\end{array}\right]\, .
\]

\s
\noindent
We obtain the following equality by computing the determinant of $M$ in two ways

\begin{equation}\label{sPair}
d\left|\begin{array}{cc}
a & f\\
e & g\end{array}\right|-g\left|\begin{array}{cc}
d & b\\
c & a\end{array}\right|=b\left|\begin{array}{cc}
c & f\\
G & g\end{array}\right|-f\left|\begin{array}{cc}
d & b\\
G & e\end{array}\right|.\end{equation}

\bs

\noindent
Notice that the left hand side of equation (\ref{sPair}) is our $S$-polynomial
$S$($ag-ef,ad-bc)$. If by a suitable choice of $G \in R$, the polynomials
$cg-fG$ and $de-bG$ are in $ \mathcal{D}$ and each term on the right hand side has order strictly less than $adg$, then the 
$S$-pair $S(ag-ef,ad-bc)$
has remainder zero modulo $ \mathcal{D}$ and in particular it reduces to zero modulo $ \D$ (see for instance  \cite[Section 2.9, Lemma 2]{CLO}). The choice of the polynomial $G$ and the
matrix $M$ depends on the polynomials we start with as we will explain
later.}
\end{Strategy}

\ms

\noindent
{\bf Proof of Lemma~\ref{GbB}} 
\noindent
Because of Remark \ref{C=B} it is enough to prove the statement for $D_{\a}$ any submatrix of $B_{\a'}$ with the same number of rows. 
Let $ \mathcal{D}$ denote the set of $2\times2$ minors of $D_{\mathfrak{a}}$.
We use the strategy described above to show that the
$S$-polynomials of any two elements of $ \D$ reduces to zero module $ \D$. 
To simplify notation, we will use $(l,\underline{j},s)$ to represent the entry $T_{{l},\uj, s}$ in the matrix $D_{\a}$. The leading terms of two elements of $ \D$ 
\[
h_1=\left|\begin{array}{cc}
(l_{1},\underline{i},s) & (l_{2},\underline{j},s)\\
(l_{1},\underline{i},t) & (l_{2},\underline{j},t)
\end{array}\right|  \qquad \text{and }
\qquad 
h_2=\left|\begin{array}{cc}
(l_{3},\underline{p},u) & (l_{4},\underline{q},u)\\
(l_{3},\underline{p},v) & (l_{4},\underline{q},v)
\end{array}\right|,
\]
 where $1\leq s<t\leq n$, $(l_{2},\underline{j},s)>(l_{1},\underline{i},s)$,  $1\leq u<v\leq n$, and $(l_{4},\underline{q},u)>(l_{3},\underline{p},u)$ are  $(l_{1},\underline{i},s) (l_{2},\underline{j},t)$ and $(l_{3},\underline{p},u)(l_{4},\underline{q},v)$. Notice that  $T_{{l},\uj, s} > T_{{l},\uj, t}$ if $ s<t$. These   will 
be relatively prime unless $(l_{1},\underline{i},s) =(l_{3},\underline{p},u)$, or $(l_{1},\underline{i},s)=(l_{4},\underline{q},v)$, or 
$ (l_{2},\underline{j},t)=(l_{3},\underline{p},u)$, or $ (l_{2},\underline{j},t)=(l_{4},\underline{q},v)$. Since $(l_{1},\underline{i},s)=(l_{4},\underline{q},v)$ and 
$ (l_{2},\underline{j},t)=(l_{3},\underline{p},u)$ are symmetric,  it suffices to consider three cases. 

\noindent
{\bf Case 1:} Set
\[
h_1=\left|\begin{array}{cc}
(l_{1},\underline{i},s)=a & (l_{2},\underline{j},s)=f\\
(l_{1},\underline{i},t)=e & (l_{2},\underline{j},t)=g
\end{array}\right|
\qquad \text{and }
\qquad 
h_2=\left|\begin{array}{cc}
(l_{1},\underline{i},s)=a & (l_{3},\underline{p},s)=c\\
(l_{1},\underline{i},u)=b & (l_{3},\underline{p},u)=d
\end{array}\right|\, ,
\]
where  we can assume $ $$s<u\leq t$ and $l_{1}\leq l_{2}\leq l_{3}$. We use equation (\ref{sPair})  and the matrix $M$ of Strategy ~\ref{strategy} with $G=(l_{3},\underline{p},t)$. The $S$-polynomial $S(h_1,h_2)$ reduces to zero modulo $ \mathcal{D}$ because $b$ and $e$ are smaller than any other entries of the two matrices defining $h_1$ and $h_2$. 

\noindent
{\bf Case 2:} Set
\[
h_1=\left|\begin{array}{cc}
(l_{1},\underline{i},s)=a & (l_{2},\underline{j},s)=e\\
(l_{1},\underline{i},t)=f & (l_{2},\underline{j},t)=g
\end{array}\right|
\qquad \text{and }
\qquad 
h_2=\left|\begin{array}{cc}
(l_{3},\underline{p},u)=d & (l_{1},\underline{i},u)=c\\
(l_{3},\underline{p},s)=b & (l_{1},\underline{i},s)=a
\end{array}\right|\, ,
\]
where $u<s<t$ and  $l_{3}\leq l_{1}\leq l_{2}$.  We use equation (\ref{sPair})  and the matrix $M$ of Strategy ~\ref{strategy} with $G=(l_{2},\underline{j},u)$. The $S$-polynomial $S(h_1,h_2)$ reduces to zero modulo $ \mathcal{D}$ because $b$ and $f$ are smaller than any other entries of the two matrices defining $h_1$ and $h_2$. 

\noindent
{\bf Case 3:} Set
\[
h_1=\left|\begin{array}{cc}
(l_{1},\underline{i},s)=d & (l_{2},\underline{j},s)=b\\
(l_{1},\underline{i},t)=c & (l_{2},\underline{j},t)=a
\end{array}\right|
\qquad \text{and }
\qquad 
h_2=\left|\begin{array}{cc}
(l_{3},\underline{p},u)=g & (l_{2},\underline{j},u)=e\\
(l_{3},\underline{p},t)=f & (l_{2},\underline{j},t)=a
\end{array}\right| \, ,
\]
where we can assume  $u\le s<t$ and  $l_{1}\leq l_{3}\leq l_{2}$. We use equation (\ref{sPair})  and the matrix $M$ of Strategy ~\ref{strategy} with $G=(l_{1},\underline{i},u)$. The $S$-polynomial $S(h_1,h_2)$ reduces to zero modulo $ \mathcal{D}$ because $c$ and $f$ are smaller than any other entries of the two matrices defining $h_1$ and $h_2$. 
\QED

\ms

\begin{Corollary}\label{prime} Adopt assumptions  ~\ref{Def1} and let $D_{\a}$ be any submatrix of $C_{\a}$ with the same number of rows.  The ideal  $I_2(D_{\a})$ is prime in  $S_{\a}$. 
\end{Corollary}
\demo Write $S=S_{\a}$. Notice that by Lemma~\ref{GbB} the variable $u \in S$ appearing in the first column and the last row of the matrix $D_{\a}$ does not divide any element in the generating set of ${\rm in}_{\tau}I_{2}(D_{\mathfrak{a}})$. Hence  $u$ is regular on $S/I_{2}(D_{\mathfrak{a}})$. 
After localizing at  $u$  the ideal $I_2(D_{\a})_{u}$ is isomorphic to an ideal generated by variables, which is a prime ideal in the ring $S_{u}$. Thus 
 $I_{2}(D_{\mathfrak{a}})$ is a prime ideal in
the ring $S$. 
\QED

\ms

\begin{Corollary}\label{generatorsinitialideals} Adopt assumptions  ~\ref{Def1}. The initial ideal $\mbox{in}_{\tau}(I_{2}(B_{\mathfrak{a}})$ is generated by the monomials $T_{{l_1},\ui, s}T_{{l_2},\uj, t}$ with $(l_1, \ui)<_{\tau} ( l_2, \uj)$ and $s<t$. 

\end{Corollary}
\demo The assertion follows from Lemma~\ref{GbB}. 
\QED
\bs

\noindent
{\bf Proof of Theorem~\ref{TH1}} 
\noindent
We first show that  for any sequence of $r$ positive integers $\a=a_1, \ldots, a_r$, the Rees algebra of $M$ is defined by the ideal of minors $I_2(C_{\a})$,
\[\R(M)=S/I_2(C_{\a}).\] 



Let $\mathcal{L}$ be the kernel of the epimorphism $\phi$ defined in (\ref{map}). Recall that $ I_2(C_{\a}) \subset \mathcal{L}$, where the first ideal is prime, according to Corollary~\ref{prime}, and the second ideal has dimension $n+r$, according to \cite[2.2]{SUV}. Hence  to show that equality holds it will be enough to prove that the  dimension of $I_2(C_{\a})$ is at most $n+r$. By  Remark \ref{C=B}, $T_{\a'}/I_2(B_{\a'})\cong S_{\a}/I_2(C_{\a})$. Hence   it will be enough to prove that the  dimension of $I_2(B_{\a'})$ is at most $n+r$ which is equivalent to show that the dimension of $I_2(B_{\a})$ is at most $n+r-1$. As ${\rm ht} \, I_{2}(B_{\mathfrak{a}})={\rm ht} \ \mbox{in}_{\tau}(I_{2}(B_{\mathfrak{a}}))$,
we can compute the dimension of  $\mbox{in}_{\tau}(I_{2}(B_{\mathfrak{a}})$. Let $U$ be the set of $r+n-1$ variables 
\[
U=\{\{T_{{l},0,\cdots,0}\}_{1\le l \le r}, T_{r,a_{r},\cdots,a_{r}},T_{r,a_{r},\ldots,a_{r},0,}\ldots,T_{r,a_{r},0,\cdots,0}\}\, .\] 
Consider the prime $\p \in {\rm Spec} (T)$ generated by all the variables of $T$ that are not in $U$. We claim that $\p$ is a  minimal prime over $\mbox{in}_{\tau}(I_{2}(B_{\mathfrak{a}})$.
It is clear that the image of $\mbox{in}_{\tau}(I_{2}(B_{\mathfrak{a}}))$ in the quotient
ring $T/\p$ is zero. To show the claim, consider the prime ideal $\p'\subset \p$  obtained by deleting   
 one variable $T_{l,\uj}$, with $\uj \in \J_{a_l}$,  from the generating set of $\p$. If $l<r$, then $T_{l, \uj}T_{r,\underline{0}}\in \mbox{in}_{\tau}(I_{2}(B_{\mathfrak{a}})\ \backslash \ \p'$. Hence  the image of $\mbox{in}_{\tau}(I_{2}(B_{\mathfrak{a}}))$ in the quotient
ring $T/\p'$ is not zero. If $l=r$, then $0<j_u \le \ldots \le j_{i} < a_r$ for some $u, i$ with $1 \le u \le i\le n-1$. Thus we can assume that 
\[T_{l,\uj}=T_{r, a_r, \ldots, a_r, j_i, \ldots, j_u, 0, \ldots,0 } \not\in \p'.\] Therefore the element $T_{l, \uj}$ appears at least in two columns of the matrix $B_{\a}$, namely the columns corresponding to the sequences \[
r, a_r, \ldots, a_r, j_i, \ldots, j_u, 1, \ldots,1\qquad \text{and} \qquad r, a_r, \ldots, a_r, j_i+1, \ldots, j_u+1, 1, \ldots,1. \] Hence   $T_{l,\uj}^2 \in \mbox{in}_{\tau}(I_{2}(B_{\mathfrak{a}})\ \backslash \ \p'$. Again   the image of $\mbox{in}_{\tau}(I_{2}(B_{\mathfrak{a}}))$ in the quotient
ring $T/\p'$ is not zero.  Hence  ${\rm dim} \, T/ I_{2}(B_{\mathfrak{a}})={\rm dim} \ T/\mbox{in}_{\tau}(I_{2}(B_{\mathfrak{a}}))={\rm dim } \, T/  \p=|U|=r+n-1$ as claimed, where the second equality follows as $  I_2(B_{\a}) $ is a prime ideal (see for instance \cite{KS}).

From the above follows that for any sequence of $r$ positive  integers $\a=a_1, \ldots, a_r$ the special fiber ring of $M$ is defined by the ideal of minors $I_2(B_{\a})$, indeed \[\mathcal{F}(M)=k\otimes_R \R(M)=k\otimes_{R}S/I_2(C_{\a})=T/I_2(B_{\a}).\] 

For any sequence $\a$ the ideals $I_2(C_{\a})$ and $I_2(B_{\a})$  have a \gb basis of quadrics according to Lemma~\ref{GbB}, hence both the Rees algebra and the special fiber ring of $M$  are Koszul domains. Normality follows because $\F(M)$ is a direct summand of $\R(M)$ which in turn is a direct summand of $R[t_1, \ldots, t_r]$. The latter claim can be easily seen once we consider $\R(M)$ as a $\NN^{r+1}-$graded $R[t_1, \ldots,t_r]$-algebra. \QED

\bs

In the rest of this section we study the divisor class group of the normal domain $A=T_{\a}/I_2(B_{\a})$ for any sequence $\a$. 

\s
\begin{Definition}\label{DEF0} {\rm Let $K$ be the $A$-ideal generated by all the variables appearing in the first row of $B_{\a}$, i.e. all the variables $T_{l,\uj}$ with $1 \le l\le r$ and $\uj \in \J_{a_l}'$. }
\end{Definition}

\ms

\begin{Theorem}\label{ClassGroup} The divisor class group $Cl(A)$ is cyclic
generated by $K$.\end{Theorem}
\demo  To compute the divisor class group we use Nagata's Theorem: If $W\subset A$ is a multiplicatively closed set, then there is an exact sequence of Abelian groups 
\[0 \lto U \lto Cl(A) \lto Cl(A_W) \lto 0\, ,\]
where $U$ is the subgroup of $Cl(A)$ generated by 
\[\{[\p] \ | \ \p \text{ a height one prime ideal with} \ \p \cap W \not=\emptyset\}\, . \]
We use the above theorem with $W=\{T^i_{r,a_{r},\cdots,a_{r}}=T_{r,\underline{a_r}}^{i}\ | \ i\in\mathbb{Z}\}$.  Notice $T_{r, \underline{a_r}} \not\in I_2(B_{\a})$. Hence  $T_{r, \underline{a_r}}$ is regular on $T/I_{2}(B_{\mathfrak{a}})$. 
After localizing at  $T_{r, \underline{a_r}}$  the ideal $I_2(B_{\a})_{T_{r, \underline{a_r}}}$ is isomorphic to an ideal generated by variables, and the ring $A_{T_{r, \underline{a_r}}}$ is a polynomial ring, hence factorial. Thus $Cl(A_W)=0$ and $Cl(A) = U$. 

Now we will show that $U$ is cyclic generated by $K$. Notice that $[K] \in U$ because $K$ is a prime ideal of height one containing $T_{r, \underline{a_r}}$. Clearly, $T_{r, \underline{a_r}}\in K$. To show that $A/K$ is a domain of dimension ${\rm dim}\, A -1$, let $R'=k[x_2, \ldots, x_n]$ be the polynomial ring over $k$ in $n-1$ variables, let $\m'$ be its homogeneous maximal ideal, and let $M'=\mathfrak{m}'^{a_{1}}\oplus\cdots\oplus\mathfrak{m}'^{a_{r}}$. The claim follows by Theorem~\ref{TH1} as $A/K \cong \F(M')$ and $\F(M')$ is a domain of dimension $n-1+r-1={\rm dim}\,  A -1.$

Let $P$ be the $A$-ideal generated by all the variables $T_{r,\uj}$ with $\uj \in \J_{a_r}$. Clearly, $P$ is prime. Indeed, if $r=1$, then $A/P=k$; if $r>1$, then, according to Theorem~\ref{TH1}, $A/P \cong \F(M_{\a'})$ with $\a'=a_1, \ldots, a_{r-1}$. Furthermore, if $r>1$, $P$ has height one as ${\rm dim} \, A/P= {\rm dim}\, \F(M_{\a'} )={\rm dim}\, A -1$.

Next we show that every prime ideal $\p$ in $A$ containing $T_{r, \underline{a_r}}$ contains either $K$ or $P$. Assume that $K \not\subset \p$, then there exists a variable $T_{l, \underline{i}} \not=T_{r, \underline{a_r}}$ with $\underline{i} \in \J'_{a_l}$ that is not in $\p$. Recall that the set of multi indices  $\J_{a_r}$ is ordered by $\tau$. For all $\underline{j}\in \J_{a_r}$ we show that $T_{r, \underline{j}}\in \p$  by descending induction on $\J_{a_r}$. The base case is trivial since $T_{r, \underline{a_r}}\in \p$.  Assume $T_{r,\underline{j}}\not=T_{r,\underline{a_r}} $.  Then there  exists a multi-index $\underline{s}  \in\J_{a_r}$ with $\underline{s} > \underline{j}$ such that the equality
\begin{equation*}\label{equality}
T_{r,\underline{j}}T_{l,\underline{i}}=T_{r, \underline{s}}T_{l, \underline{i}, t}
\end{equation*}
holds in $A$ for some integer $t\ge 2$.  Notice that $T_{l, \underline{i}, t}$ is well defined because $\underline{i} \in \J'_{a_l}$. By induction $T_{r, \underline{s}} \in \p$, thus $T_{r, \underline{j}} \in \p$ since $\p$ is prime and 
$T_{l,\underline{i}}\not\in \p$.  Hence $P \subset \p$.

If $r=1$ the above inclusion implies that ${\rm ht}\, \p >1$ and hence $U$ is cyclic generated by $K$. If $r>1$, then $U$ is generated by $P$ and $K$. We conclude by showing that $[P]=a_r[K]$, or equivalently, $P=(T_{r, \underline{a_r}}):_AK^{a_r}$. Since both ideals have height one and $P$ is prime it is enough to prove the inclusion $P\subset(T_{r, \underline{a_r}}):_AK^{a_r}$.  The latter follows from the equation
\begin{equation}\label{induction}
T_{r,\underline{j}} K ^{a_r-j_1} \in (T_{r, \underline{a_r}})\, .
\end{equation}
We prove equation (\ref{induction}) using descending induction on $j_1$ with $0\le j_1 \le a_r$.  The base case is trivial since $j_1=a_r$ implies $T_{r,\underline{j}}=T_{r, \underline{a_r}}$. If $j_1 < a_r$, then the variable $T_{r,\underline{j}}$ appears in a row $s$ of $B_{\a}$ with $s>1$. Thus for any element $\lambda \in K$ there exists $\beta \in S$ such that the equality 
\[T_{r, \underline{j}} \lambda=T_{r, j_{n-1},\cdots, j_{s+1}, j_s+1, \cdots, j_{1}+1} \beta\] holds in $A$. Now the claim follows by induction. 
\QED

\bs

According to  Theorem~\ref{ClassGroup} the classes of the divisorial ideals $K^{(\delta)}$ and $P^{(\delta)}$, the $\delta$-th symbolic power of $K$ and $P$ respectively, constitute $Cl(A)$.  In the next theorem we exhibit a monomial generating set for $K^{(\delta)}$ for $\delta \ge 1$.
As in \cite{KPU1} and \cite{S} 
we identify $K^{(\delta)}$ with
a graded piece of $A$. We put a new grading on the ring $T$,
\[
\mbox{Deg}(T_{l, \uj})=j_{1}.
\]
Notice that $I_{2}(B_{\mathfrak{a}})$ is an homogeneous ideal with respect
to this grading. Thus Deg  induces a grading on $A$. Let $A_{\geq\delta}$
be the ideal generated by all monomials $m$ in $A$ with 
${\rm Deg}(m) \ge \delta$. 
\bs
\begin{Theorem}\label{KnThm} The $\delta$-th symbolic power of $K$, $K^{(\delta)}$, equals the monomial ideal $A_{\geq\delta}$.
\end{Theorem}
\demo One proceeds as in \cite[1.5]{KPU2}.
\QED

\bs
\bs
\section{The Blowup algebras of truncations of complete intersections}

\bs

In this section our goal is to compute explicitly the defining equations of the blowup
algebras of truncations of  complete intersections.

\begin{Assumptions}\label{DEF}{\rm 
Let $R=k[x_1,\ldots, x_n]$ be a polynomial ring over a field $k$ with homogeneous maximal ideal ${\mathfrak m}$.  Let $r$ be an integer with $1\le r \le n$. 
Let $f_{1}$, \ldots, $f_{r}$ be a homogeneous regular sequence
in $R$ of degree $d_{1}\ge \ldots \ge d_{r}$. Let $d\ge d_1$ be an integer, write  $a_i=d-d_i$.  Let $I$ be the truncation of the complete intersection $(f_1,\ldots, f_r)$ in degree $d$, i.e. the $R$-ideal generated by $\{ \mathfrak{m}^{a_{i}}f_{i} \ | \ 1\le i \le r\}$ ,
\[
I=(f_1, \ldots, f_r)_{\ge d}=(\mathfrak{m}^{a_{1}}f_{1},\ldots, \mathfrak{m}^{a_{r}}f_{r}).
\]}
\end{Assumptions}

\bs

\begin{Theorem}\label{RI} Adopt assumptions  ~\ref{DEF} and write $M=\m^{a_1}\oplus \ldots \oplus \m^{a_r}$. There is a  short exact sequence 
\[
0\rightarrow Q\rightarrow\mathcal{R}(M)\rightarrow\mathcal{R}(I)\rightarrow0
\] 
where $Q$ is a prime ideal of height $r-1$ in the normal domain $\R(M)$. 

\end{Theorem}
\demo
The natural map 
\begin{eqnarray*}
M \qquad & \stackrel{\zeta}\lto & \  I=\mathfrak{m}^{a_{1}}f_{1}+\ldots +\mathfrak{m}^{a_{r}}f_{r}\rightarrow0\\
(u_1,\ldots, u_r) & \mapsto & \qquad  u_1f_{1}+ \ldots + u_rf_{r}
\end{eqnarray*}
 induces a surjection on the level of Rees algebras
\[
\Psi: \mathcal{R}(M)\lto\mathcal{R}(I) 
\]
and ${\rm Ker} \ \Psi$  is a prime of height $r-1$ since $\R(M)$ and $\mathcal{R}(I)$
are domains and $\dim\mathcal{R}(M)=n+r=\dim\mathcal{R}(I)+r-1$.
\QED
\bs

\begin{Definition}\label{kernel} {\rm Adopt assumptions  ~\ref{DEF}. Let  $\a$ be the sequence $a_1, \ldots, a_r$ and let $S$ be the polynomial ring $S_{\a}$. 
Consider the algebra epimorphism $\chi$ obtained by the composition of the two algebra epimorphisms
\[ \phi: S \lto \R(M) \qquad \text{and} \qquad  \Psi:\R(M) \lto \R(I)
\]
with $\phi$ as in (\ref{map}) and let $\mathcal{A}$ be the $S$-ideal defined by the short
exact sequence 
\[
0\rightarrow\mathcal{A}\lto S  \stackrel{\chi}\lto\mathcal{R}(I)\rightarrow0\, .
\]}
\end{Definition}

\bs

Let $\tau$ be the order on $T_{l,\uj}$ defined in Section \ref{section1}. Notice that, through the algebra epimorphism $\chi$, the order $\tau$ induces an order on the set \[\C=\left\{\ux^{a_l,\uj}f_{l} \ |\  1 \le l \le r, \uj\in \J_{a_l} \right\}\] of  generators of $I$. Let $\varphi$ be the 
presentation matrix of $I$ over $R$ with respect to $\C$. In the following remark we give a resolution of $I$ when $r=2$. In particular, we describe explicitly $\varphi$. 

\begin{Remark}\label{varphi}{\rm To compute a resolution of $I=(f_1,f_2)_{\ge d}$ we first truncate the Koszul complex of  $f_1, f_2$ in degree $d$. 
\begin{enumerate}
\item If  $d_2 \le d\le d_1+d_2$ we obtain the short exact sequence:
\[0 \lto R(-d_1-d_2+d) \stackrel{\rho}\lto M \stackrel{\zeta}\lto I \lto 0 \, ,
\]
where $\zeta$ is the natural surjection described in Theorem~\ref{RI} and $\rho= \left[
\begin{array}{cc}
 -f_2  \\ 

\  \ f_1  \\
\end{array}
\right]$. Set $\varphi''$ to be the map that rewrites  
$ \left[
\begin{array}{cc}
 -f_2  \\ 

\  \ f_1  \\
\end{array}
\right]$ in terms of the $k$-basis $\B$ of $M$ ordered with the order induced by $\tau$.  The Eagon-Northcott complex gives us an $R$-resolution $F_{\bullet}$ of $M$, while the complex $G_{\bullet}$ that is trivial everywhere except in degree $0$ gives us an $R$-resolution of $R(-d_1-d_2+d)$. The map $\rho$ can be trivially lifted to a morphism of complexes 
\[\rho_{\bullet}: G_{\bullet} \lto F_{\bullet}\] that is trivial in positive degree and is $\varphi''$ in degree zero.
\item  If  $ d\ge d_1+d_2+1$ we obtain the short exact sequence:
\[0 \lto \m^{d-d_1-d_2} \stackrel{\rho}\lto M \stackrel{\zeta}\lto I \lto 0 \, ,
\]
where $\zeta$ is the natural surjection described in Theorem~\ref{RI} and $\rho= \left[
\begin{array}{cc}
 -f_2  \\ 

\  \ f_1  \\
\end{array}
\right]$. Write $\E$ for the $k$-bases of $\m^{d-d_1-d_2}$  ordered with the order induced by $\tau$. Set $\varphi''$ to be  the map that rewrites  
$ \E \left[
\begin{array}{cc}
-f_2  \\ 

\ \ f_1  \\
\end{array}
\right]$  
in terms of the $k$-basis $\B$ of $M$ ordered with the order induced by $\tau$. From the Eagon-Northcott complex we obtain $R$-resolutions $F_{\bullet}$ and $G_{\bullet}$ of $M$ and $\m^{d-d_1-d_2}$, respectively. The map $\rho$ can be  lifted to a morphism of complexes 
\[\rho_{\bullet}: G_{\bullet} \lto F_{\bullet}\] with $\varphi''$ in degree zero.  
 \end{enumerate}
In both cases the mapping cone $C(\rho_{\bullet})$ is a non-minimal free resolution of $I$. In particular, the presentation matrix of $I$ with respect to $\C$ is $\varphi=[\varphi',\varphi'']$ where $\varphi'$ is the matrix presenting $M$ with respect to $\B$. 

 }
\end{Remark}

 In the following remark we describe explicitly $\varphi$ for any $r$ when $d\ge d_1+d_2$. 

\begin{Remark}\label{h2} {\rm  Adopt assumptions  ~\ref{DEF}. To compute the presentation matrix for  $I=(f_1, \ldots, f_r)_{\ge d}$ we  need to truncate the Koszul complex $K_{\bullet}(f_1, \ldots, f_r)=\bigwedge Re_1 \oplus \ldots \oplus Re_r$. If $d \ge d_1 +d_2$ then $d\ge d_i+d_j$ for any $1\le i < j\le r$ and we have:

\[
\ldots \lto \bigoplus_{1\le i<j\le r}\m^{d -d_i-d_j} e_i \wedge e_j \lto M=\bigoplus_{i=1}^{r} \m^{a_i} \lto I \lto 0\, .
\]
Write $\B$ and $\E_{i,j}$ for the $k$-bases of $M$ and of $\m^{d-d_i-d_j}$, respectively,   ordered with the order induced by $\tau$. Set $\varphi''$ to be  the map that rewrites  
$ \E_{i,j} (-f_je_i+f_ie_j)$   in terms of $\B$ for all $1\le i <j\le r$. The presentation matrix of $I$ with respect to $\C$ is $\varphi=[\varphi',\varphi'']$ where $\varphi'$ is the matrix presenting $M$ with respect to $\B$.
}
\end{Remark}
\ms

\begin{Definition}\label{h}{\rm Adopt assumptions  ~\ref{DEF}. Let 
$h_1, \ldots, h_t \in S$  be the homogeneous  polynomials obtained by the matrix multiplication
$[\underline{T}]\varphi''$, where $\varphi''$ is the matrix described in Remark \ref{varphi} and \ref{h2}. Think of $S$ as
a naturally bigraded ring with ${\rm deg} \ x_i=(1,0)$ and  ${\rm deg} \ T_j=(0,1)$. If  $r=2$ and $d_1 \le d\le d_1+d_2$, then $t=1$ and $h_1$ has bidegree $(\delta,1)$ with $\delta=d_2-a_1=d_1-a_2=d_1+d_2-d \ge 0$. For any $r$, if $ d\ge d_1+d_2$ then $t=
\displaystyle\sum_{1\le i<j\le r}{\sigma_{i,j}+n-1\choose n-1}$ with $\sigma_{i,j}=d-d_i-d_j\ge 0$ and the $h_k$'s have bidegree $(0,1)$.}
\end{Definition}

\ms

\begin{Proposition}\label{hprime} Adopt assumptions  ~\ref{DEF}, \ref{h} with $ d\ge d_1+d_2$ then 
$ (h_1, \ldots, h_t) \R(M)$ is a non-zero prime ideal of height $\ge r-1$. 
\end{Proposition}
\demo  Let $1\le i < j \le r$. Write \[f_i=\sum_{\underline{k} \in \J_{d_i}} \lambda_{\underline{k}} \ux^{d_i, \uk}  \qquad \text{and} \qquad f_j=\sum_{\underline{k} \in \J_{d_j}} \alpha_{\underline{k}} \ux^{d_j, \uk} \]  Let $\sigma_{i,j}= d-d_i-d_j\ge 0$. As in Definition~\ref{N0}, denote the elements of the basis 
 $\E_{i,j}$ of $\m^{\sigma_{i,j}}$ with $\ux^{\sigma_{i,j},\ut,s} $, $\ut \in \J_{\sigma_{i,j}}'$. Since we have a one to one correspondence between the $h_k$'s and the  elements of $\E_{i,j}$, we write $h_k$ as $h_{\sigma_{i,j}, \ut, s}$. Let $\H$ be the $S$-ideal generated by $\{h_{\sigma_{i,j}, \ut, s} \ | \  1\le i < j \le r, \ut \in \J'_{\sigma_{i,j}}, 1\le s\le n\}$. We obtain 
 \[ \ux^{\sigma_{i,j},\ut,s} f_i= \sum_{\underline{k} \in \J_{d_i}} \lambda_{\underline{k}} \ux^{a_i, \uk+\ut,s} \qquad \text{and} \qquad  \ux^{\sigma_{i,j},\ut,s} f_j= \sum_{\underline{k} \in \J_{d_j}} \alpha_{\underline{k}} \ux^{a_j, \uk+\ut,s}\, , \]
hence
\[h_{\sigma_{i,j}, \ut, s}=  \sum_{\underline{k} \in \J_{d_i}} \lambda_{\underline{k}} T_{i, \uk+\ut,s} - \sum_{\underline{k} \in \J_{d_j}} \alpha_{\underline{k}} T_{j, \uk+\ut,s}\, .
\]
Let $\H_j$ be the $S$-ideal generated by $\{h_{\sigma_{1,l}, \ut, s} \ | \  2\le l\le  j,  \ \ut \in \J'_{\sigma_{1,l}}, \ 1\le s\le n\}$. 
We show by induction on $j$ with $2\le j \le r$ that the ideal $\H_j\R(M)$ is prime of height $\ge j-1$. Let  $j=2$. Notice that for each $\ut \in \J'_{\sigma_{1,2}}$,  the column $[h_{\sigma_{1,2}, \ut, 1}, \ldots, h_{\sigma_{1,2}, \ut, n}]^{tr}$ is a linear combination of columns of $C_{\a}$. Write $E_{\a}$ for the  $n \times (1+\sum_{q=1}^{r}\binom{a_{q}+n-2}{n-1})$ matrix obtained by $C_{\a}$ by substituting 
 ${\sigma_{1,2} +n-2\choose n-1}$ columns  with  $[h_{\sigma_{1,2}, \ut, 1}, \ldots, h_{\sigma_{1,2}, \ut, n}]^{tr}$,  $\ \ut \in \J'_{\sigma_{1,2}}$. The  two $S$-ideals $I_2(C_{\a})+ \H_2$ and $I_2(E_{\a})+\H_2$ are equal. The ideal $\H_2\R(M)$ is prime according to  Corollary~\ref{prime} as 
 \[S/(I_2(E_{\a})+\H_2) \cong T'/I_2(D_{\a})\]
 for some polynomial ring $T'$ and $D_{\a}$ a suitable submatrix of $C_{a}$ with $n$ rows. Now degree considerations show that $h_{\sigma_{1,2}, \ut, s} \not\in I_2(C_{\a})$, hence $\H_2\R(M)\not=0$.  Thus the ideal $\H_2\R(M)$ is prime of height at least one. 
 
 Let $2 \le l \le  j$ and assume by induction that $\H_{j-1}\R(M)$ is prime of height $\ge j-2$, we show $\H_{j}\R(M)$ is prime of height $\ge j-1$. For each $\ut \in \J'_{\sigma_{1,l}}$,  the column $[h_{\sigma_{1,l}, \ut, 1}, \ldots, h_{\sigma_{1,l}, \ut, n}]^{tr}$ is a linear combination of columns of $C_{\a}$. Write $E_{\a}$ for the  $n \times (1+\sum_{q=1}^{r}\binom{a_{q}+n-2}{n-1})$ matrix obtained by $C_{\a}$ by substituting 
 $\sum_{l=2}^{j}{\sigma_{1,l} +n-2\choose n-1}$ columns  with  $[h_{\sigma_{1,l}, \ut, 1}, \ldots, h_{\sigma_{1,l}, \ut, n}]^{tr}$,  $2\le l \le j$ and $\ \ut \in \J'_{\sigma_{1,l}}$. The  two $S$-ideals $I_2(C_{\a})+ \H_l$ and $I_2(E_{\a})+\H_l$ are equal. The ideal $\H_j\R(M)$ is prime according to  Corollary~\ref{prime} as 
 \[S/(I_j(E_{\a})+\H_j) \cong T'/I_2(D_{\a})\]
 for some polynomial ring $T'$ and $D_{\a}$ a suitable submatrix of $C_{a}$ with $n$ rows.  Notice that for each $\ut \in \J'_{\sigma_{1,j}}$,  the column $[h_{\sigma_{1,j}, \ut, 1}, \ldots, h_{\sigma_{1,j}, \ut, n}]^{tr}$ is a linear combination of  a subset of the columns $\{(1,\ut) \ | \ \ut \in  \J'_{{a_1}}\}$ and $\{(j,\ut) \ | \ \ut \in  \J'_{{a_j}}\}$ of $C_{\a}$, while the column $[h_{\sigma_{1,l}, \ut, 1}, \ldots, h_{\sigma_{1,l}, \ut, n}]^{tr}$ is a linear combination of  a subset of the columns $\{(1,\ut) \ | \ \ut \in  \J'_{{a_1}}\}$ and $\{(l,\ut) \ | \ \ut \in  \J'_{{a_l}}\}$ of $C_{\a}$ with $2\le l \le j-1$. Hence degree considerations show that $h_{\sigma_{1,j}, \ut, s} \not\in (\H_{j-1},I_2(C_{\a}))$, hence $\H_j\R(M)\not=0$.  Thus the ideal $\H_j\R(M)$ is prime of height $\ge j-1$. 
 
Using the same argument  one can show that $\H\R(M)$ is a prime ideal and its height is at least $r-1$ as $\H\R(M) \supset \H_r\R(M)$ . 
  \QED  
\ms

\begin{Remark}{\rm  If $d\ge d_1+d_2-1$, then the ideal $I=(f_1, \ldots, f_r)_{\ge d}$ has a linear resolution.
The Rees algebra of linearly presented ideals of height 2 has been described explicitly in terms of generators and relations in \cite{MU} under the additional assumption that $I$ is perfect. However, if  $r<n$, the truncations of codimension $r$ complete intersections are never perfect. For large $d$, the Rees algebra  $\R(I)$ is Cohen-Macaulay as shown in \cite{HT}. But the defining equations of $\R(I)$ were unknown, we give them explicitly in Theorem~\ref{anys}. If  $r=2$ we prove that for $d\ge d_1+d_2-1$ the Rees ring $\R(I)$ is a Koszul domain (see Corollary~\ref{koszul}). 
If $r\ge 3$, we prove that  $\R(I)$ is a Koszul domain for $d\ge d_1+d_2$ (see Corollary~\ref{koszul}).
In addition in \cite{LP} we study the depth and regularity of  the blowup algebras of $I$. }
\end{Remark}

\ms

\begin{Theorem}\label{anys} Adopt assumptions  ~\ref{DEF} and  \ref{h} with $d\ge d_1+d_2$ then \[\R(I)=\R(M_{\a})/(h_1, \ldots, h_t)= S_{\a}/(I_2(C_{\a}), h_1, \ldots, h_t) \, .\]
\end{Theorem}
\demo Adopt the notation of the proof of Proposition \ref{hprime}.   Let $\H$ be the $S$-ideal generated by $\{h_{\sigma_{i,j}, \ut, s} \ | \ 1 \le i < j \le r, \ut \in \J'_{\sigma}, 1\le s\le n\}$. According to Remark~\ref{h2},  we have $\H\subset \A$.  Hence $\H\R(M) \subset \A\R(M)$ where the first ideal is a non-zero prime ideal  of height $\ge r-1$ by Proposition \ref{hprime} and the second one has height $r-1$. 
\QED

\bs

Assume $r=2$. In the following theorem we express the Rees algebra of $I$ as defined by a divisor on the Rees algebra of the module $M=\m^{a_1}\oplus \m^{a_2}$ that we computed  explicitly in the previous section. Indeed, for $r=2$ the prime ideal $Q \in {\rm Spec} (\R(M))$ of Theorem \ref{RI}  gives rise to an element of the divisor class   $Cl(\R(M))$. This group has been studied explicitly in Theorem~\ref{ClassGroup}: it is cyclic generated by the prime ideal $L$, where $L$ be the $\R(M)$-ideal generated by all the variables appearing in the first row of $C_{\a}$. In the next theorem we identify    
for which $s$ the ideal $L^{(s)}$ is isomorphic to $Q$.
\s
 
\begin{Definition}\label{DEFL} {\rm Let $L$ be the $\R(M)$-ideal generated by all the variables appearing in the first row of $C_{\a}$. }
\end{Definition}

\s
 
We will use the convention that the (symbolic) power of any element or ideal with nonpositive
exponent is one or the unit ideal, respectively.
 \ms
 
\begin{Theorem}\label{main} Adopt assumptions  ~\ref{DEF}, \ref{h}, and \ref{DEFL} with $r=2$ then
\[x_1^{\delta}\A\R(M)=(h_1, \ldots, h_t)L^{(\delta)}\, .
\]
In particular,
\begin{enumerate}
\item
If $d_1 \le d \le d_1+d_2-1$ then the $\R(M)$-ideals $x_1^{\delta}\A\R(M)$ and $h_1L^{(\delta)}$ are equal and the bigraded $\R(M)$-modules $\A$ and $L^{(\delta)}(0,-1)$ are isomorphic.
\item If $d\ge d_1+d_2$ then \[\R(I)=S_{\a}/(I_2(C_{\a}), h_1, \ldots, h_t)\, .\]
\end{enumerate}
\end{Theorem}
\demo Notice that the first statement follows from (1) and (2) and (2) has been proven in Theorem \ref{anys}. Thus it will be enough to prove (1).
Write $h=h_1$. One proceeds as in \cite[1.11]{KPU2}. For clarity we rewrite part of the proof here since there are some minor differences and in \cite[1.11]{KPU2} the integer  $\delta$ was assumed to be  $ \ge 2$. Degree considerations show that $x_1$ is not in $I_2(C_{\a})$. The ideal $I_2(C_{\a})$ is prime, so $x_1^{\delta}$ is also not in $I_2(C_{\a})$. The second assertion in (1) follows from the first as $x_1^{\delta}$ has bi-degree $(\delta,0)$ and $h_1$ has bi-degree $(\delta,1)$. Write $\overline{\phantom{x}}$ to mean image in $\R(M)$. We prove the equality  $x_1^{\delta}\A\R(M)=hL^{(\delta)}$ by showing that $\overline{\A}=({\bar h}/{\bar x_1^{\delta}})L^{(\delta)}$, where the fraction is taken in  the quotient field $Q$ of $\R(M)$.

Notice that $(\bar x_1^i):_QL^{(i)}=({\bar x_1},\ldots,  {\bar x_n})^i $. This follows as in the proof of claim (1.12) in \cite[1.11]{KPU2}. Furthermore $\bar h \in \overline{\m}^{\delta}=  (\bar x_1^{\delta}):_QL^{({\delta})}$. Thus, $\bar hL^{(\delta)}\subseteq \bar x_1^{\delta} \R(M)$. Define $D$ to be the ideal   $(\bar h/\bar x_1^{\delta})L^{(\delta)}$ of $\R(M)$.   At this point, we see that the ideal $D$ is either zero or divisorial.

To show that $D$ is not zero and to establish
the equality $\overline{\A}=D$, it suffices to prove that $\overline{\A}\subseteq D$, because $\overline{\A} $ is a height one  prime ideal of $\R(M)$. This is the only part where the argument differs from  \cite[1.11]{KPU2}. Notice that $\bar h\in D$ as $\bar x_1\in L$. For every $w\in \m$, one has $I_w=(f_1,f_2)_w$. Therefore, $R[It]_w=R[(f_1,f_2)t]_w$ and we obtain $(\bar h)_w=\overline{\A}_w$. It follows that $\bar h\neq 0$  and   $\overline{\A}_w\subseteq D_w$. The rest of the proof follows as in \cite[1.11]{KPU2}.
 \QED
 
 \ms
 \begin{Corollary}\label{fiber}Adopt assumptions  ~\ref{DEF}, \ref{h}, and \ref{DEF0}. 
 \begin{itemize}
 \item[(a)]  If $d\ge d_1+d_2$ then \[\F(I)=\F(M_{\a})/(h_1, \ldots, h_t)= T_{\a}/(I_2(B_{\a}), h_1, \ldots, h_t) \, .\]

 \item[(b)] If $r=2$  and $d_1 \le d \le d_1+d_2-1$ then
\[\F(I)=T_{\a}/{\mathcal K}\, \qquad \mbox{with} \qquad {\mathcal K} \cong  K^{(\delta)}(-1) .
\]

 \end{itemize}
 \end{Corollary}
 \demo The proof follows from Theorem \ref{anys} and Theorem \ref{main} and the fact that $\F(I)= k \otimes \R(I)$.  Also the same argument  as in \cite[4.2]{KPU2} shows the isomorphism ${\mathcal K} \cong  K^{(\delta)}(-1)$.
 \QED

\ms

\begin{Corollary}\label{koszul}  Adopt assumptions  ~\ref{DEF}.  If  $d \ge d_1+d_2$ or if $r=2$ and $d=d_1+d_2-1$, then $\R(I)$ and $\F(I)$ are  Koszul algebras.
\end{Corollary}
\demo  We will prove the statement for the Rees ring. The same proof works for the special fiber ring using Corollary \ref{fiber}. 
If $d\ge d_1+d_2$, then  $\R(I)=S/(I_2(C_{\a}), h_1, \ldots, h_t)$ according to Theorem~\ref{anys}, and according to the proof of Corollary~\ref{prime} the latter ring is isomorphic to $T'/I_2(D_{\a})$
 for some polynomial ring $T'$ and $D_{\a}$ a suitable submatrix of $C_{a}$ with $n$ rows. But this ring is Koszul as it has a \gb of quadrics by Lemma \ref{GbB}.

If $r=2$ and  $d=d_1+d_2-1$, then $\delta=1$ and $\R(I)= \R(M)/{\mathcal A}\R(M)\cong S/(I_2(C_{\a}), L)$ according to Theorem~\ref{main}(1). Let $R'=k[x_2, \ldots, x_n]$ be the polynomial ring over $k$ in $n-1$ variables, let $\m'$ be its homogenoeus maximal ideal, and let $M'=\mathfrak{m}'^{a_{1}}\oplus\mathfrak{m}'^{a_{2}}$. The last assertion now follows from Theorem~\ref{TH1} as \[\R(I)= \R(M)/{\mathcal A}\R(M)\cong S/(I_2(C_{\a}), L)\cong \R(M')\, .\]
\QED

\bs

\end{document}